\documentclass[leqno, 12pt]{article}
\usepackage{amsmath,amsfonts,amsthm,amssymb,indentfirst}

\newtheorem{theorem}{Theorem}
\newtheorem{lemma}[theorem]{Lemma}
\newtheorem{definition}[theorem]{Definition}
\newtheorem{corollary}[theorem]{Corollary}
\newtheorem{proposition}[theorem]{Proposition}

\begin{document}

\title{Commutator rings}
\author{Zachary Mesyan}

\maketitle

\begin{abstract}
A ring is called a commutator ring if every element is a sum of
additive commutators. In this note we give examples of such rings.
In particular, we show that given any ring $R$, a right $R$-module
$N$, and a nonempty set $\Omega$, $\mathrm{End}_R
(\bigoplus_{\Omega} N)$ and $\mathrm{End}_R (\prod_{\Omega} N)$
are commutator rings if and only if either $\Omega$ is infinite or
$\mathrm{End}_R (N)$ is itself a commutator ring. We also prove
that over any ring, a matrix having trace zero can be expressed as
a sum of two commutators.
\end{abstract}

\section{Introduction}

In 1956 Irving Kaplansky proposed twelve problems in ring theory
(cf.~\cite{Kaplansky}). One of these was whether there is a
division ring in which every element is a sum of additive
commutators. This was answered in the affirmative two years later
by Bruno Harris in~\cite{Harris}. However, it seems that to this
day not much is known about rings in which every element is a sum
of commutators, which we will call commutator rings.

The purpose of the first half of this note is to collect examples
of such rings. For instance, we will show that given any ring $R$,
a right $R$-module $N$, and a nonempty set $\Omega$,
$\mathrm{End}_R (\bigoplus_{\Omega} N)$ and $\mathrm{End}_R
(\prod_{\Omega} N)$ are commutator rings if and only if either
$\Omega$ is infinite or $\mathrm{End}_R (N)$ is itself a
commutator ring. We will also note that if $R$ is a nonzero
commutative ring and $n$ is a positive integer, then the Weyl
algebra $A_n(R)$ is a commutator ring if and only if $R$ is a
$\mathbb Q$-algebra.

Along the way, we will give an alternate characterization of right
$R$-modules $M$ that can be written in the form
$\bigoplus_{i\in\Omega} N_i$, where $\Omega$ is an infinite set,
and the right $R$-modules $N_i$ are all isomorphic. Specifically,
$M$ is of this form if and only if there exist $x, z \in
\mathrm{End}_R(M)$ such that $zx = 1$ and $\bigcup_{i = 1}^\infty
\mathrm{ker}(z^i) = M$.

The last section of this note is devoted to commutators in matrix
rings. In \cite{A&M} Albert and Muckenhoupt showed that if $F$ is
a field and $n$ is a positive integer, then every matrix in
$\mathbb{M}_n (F)$ having trace zero can be expressed as a
commutator in that ring.  (This was first proved for fields of
characteristic zero by Shoda in \cite{Shoda}.) The question of
whether this result could be generalized to arbitrary rings
remained open for a number of years, until M. Rosset and S. Rosset
showed in~\cite{Rosset2} that it could not. (An example
demonstrating this will also be given below.) However, we will
prove that every matrix having trace zero can be expressed as a
sum of two commutators, generalizing a result of M. Rosset
in~\cite{Rosset} (unpublished) for matrices over commutative
rings. As a corollary, we also obtain a generalization to
arbitrary rings of the result in~\cite{Harris} that a matrix over
a division ring is a sum of commutators if and only if its trace
is a sum of commutators. On a related note, in~\cite{A&R} Amitsur
and Rowen showed that if $R$ is a finite-dimensional central
simple algebra over a field $F$, then every element $r \in R$ such
that $r \otimes 1$ has trace zero in $R \otimes_F \bar{F} \cong
\mathbb{M}_n (\bar{F})$ is a sum of two commutators, where
$\bar{F}$ is the algebraic closure of $F$. (See
also~\cite{Rosset1}.)

I am grateful to George Bergman for his numerous comments and
suggestions on earlier drafts of this note. Also, many thanks to
Lance Small for his comments and for pointing me to related
literature.

\section{Definitions and examples}

Given a unital associative ring $R$ and two elements $x,y \in R$,
let $[x, y]$ denote the commutator $xy - yx$. We note that $[x,
y]$ is additive in either variable and $[x, yz] = [x, y]z + y[x,
z]$, $[zx, y]  = [z, y]x + z[x, y]$ are satisfied by all $x,y,z
\in R$ (i.e., $[x, - ]$ and $[ -, y]$ are derivations on $R$). Let
$[R, R]$ denote the additive subgroup of $R$ generated by the
commutators in $R$, and given an element $x \in R$, let $[x, R] =
\{[x, y] : y \in R \}$. If $n$ is a positive integer, we will
denote the set of sums of $n$ commutators in elements of $R$ by
$[R, R]_n$. For convenience, we define $[R, R]_0 = \{0\}$.
Finally, right module endomorphisms will be written on the left of
their arguments.

\begin{definition}
A ring $R$ is called a \emph{commutator ring} if $R = [R, R]$.
\end{definition}

In~\cite{Harris} and~\cite{Lazerson} examples of commutator
division rings are constructed. Also, it is easy to see that
finite direct products and homomorphic images of commutator rings
are commutator rings.

\begin{proposition}\label{superrings}
If $R \subseteq S$ are rings such that $R$ is a commutator ring
and $S$ is generated over $R$ by elements centralizing $R$, then
$S$ is also a commutator ring.
\end{proposition}

\begin{proof}
Given an element $a \in S$, we can write $a = \sum_{i = 1}^m r_i
s_i$, where $r_i \in R$, the $s_i \in S$ centralize $R$, and $m$
is some positive integer. Since $R$ is a commutator ring, for each
$r_i$ there are elements $y_{ij}, x_{ij} \in R$ and a positive
integer $m_i$ such that $r_i = \sum_{j=1}^{m_i} [x_{ij}, y_{ij}]$.
Then
$$a = \sum_{i = 1}^m r_i s_i = \sum_{i = 1}^m
(\sum_{j=1}^{m_i} [x_{ij}, y_{ij}])s_i = \sum_{i = 1}^m
(\sum_{j=1}^{m_i} [x_{ij}, y_{ij}s_i]).$$
\end{proof}

This proposition implies that, for example, matrix rings, group
rings, and polynomial rings over commutator rings are commutator
rings.  Furthermore, given a commutative ring $K$ and two
$K$-algebras $R$ and $S$, such that $R$ is a commutator ring, $R
\otimes_K S$ is again a commutator ring.

Given a ring $R$, a set of variables $X$, and a set of relations
$Y$, we will denote the $R$-ring presented by $R$-centralizing
generators $X$ and relations $Y$ by $R \langle X : Y \rangle$.

\begin{definition}
Let $R$ be a ring.  Then $A_1(R) = R \langle x, y : [x, y] = 1
\rangle$ is called the \emph{(first) Weyl algebra} over $R$.
Higher Weyl algebras are defined inductively  by $A_n(R) =
A_1(A_{n-1}(R))$.  More generally, given a set $I$, we will denote
the ring
$$R \langle \{x_i\}_{i \in I}, \{y_i\}_{i \in I} : [x_i,
x_j] = [y_i, y_j] = [x_i, y_j] = 0 \ \mathrm{for} \ i \neq j, \
\mathrm{and} \ [x_i, y_i] = 1 \rangle$$
by $A_I(R)$.
\end{definition}

\begin{proposition}\label{A_I(R)}
For any ring $R$ and any infinite set $I$, $A_I(R) = [A_I(R),
A_I(R)]_1$. In particular, any ring can be embedded in a
commutator ring.
\end{proposition}

\begin{proof}
Let $s \in A_I(R)$ be any element.  When written as a polynomial
in $\{x_i\}_{i \in I}$ and $\{y_i\}_{i \in I}$, $s$ contains only
finitely many variables $y_i$, so there exists some $n \in I$ such
that $y_n$ does not occur in $s$. Then $[x_n,\ s] = 0$, and hence
$s = 1 \cdot s = [x_n, y_n]s = [x_n, y_ns] - y_n[x_n, s] = [x_n,
y_ns]$.
\end{proof}

Harris used this construction in~\cite{Harris} to produce a
commutator division ring. Specifically, he took $R$ to be a field
and then used essentially the method above to show that the
division ring of fractions of $A_I(R)$ is a commutator ring.

Before discussing the case when $I$ is a finite set, we record two
well known facts that will be useful.

\begin{lemma}\label{matrix trace}
Let $R$ be any ring and $n$ a positive integer. If $A, B \in
\mathbb{M}_n (R)$, then $\mathrm{trace}([A, B]) \in [R, R]_{n^2}$.
Conversely, given any $r \in [R, R]_{n^2}$, there exist matrices
$A, B \in \mathbb{M}_n (R)$ such that $\mathrm{trace}([A, B]) =
r$.
\end{lemma}

\begin{proof}
Let $A, B \in \mathbb{M}_n (R)$, and write $A = (a_{ij})$ and $B =
(b_{ij})$. Then $$\mathrm{trace} (AB - BA) = \sum_{i=1}^n
\sum_{k=1}^n a_{ik} b_{ki} - \sum_{i=1}^n \sum_{k=1}^n b_{ik}
a_{ki}.$$ Now, for all $i$ and $j$ with $1 \leq i, j \leq n$,
$a_{ij} b_{ji}$ appears exactly once in the first term of this
expression, and $b_{ji} a_{ij}$ appears exactly once in the second
term. Hence, $$\mathrm{trace} (AB - BA) = \sum_{i=1}^n
\sum_{j=1}^n [a_{ij}, b_{ji}].$$

For the second statement, given $r \in [R, R]_{n^2}$, we can write
$r = \sum_{i=1}^n \sum_{j=1}^n [a_{ij}, b_{ji}]$ for some $a_{ij},
b_{ji} \in R$, after relabeling. Setting $A = (a_{ij})$ and $B =
(b_{ij})$, we have $\mathrm{trace}([A, B]) = r$, by the previous
paragraph.
\end{proof}

\begin{proposition}\label{fd algerbra}
Let $F$ be a field and $R$ a finite-dimensional $F$-algebra. Then
$R \neq [R, R]$.
\end{proposition}

\begin{proof}
Suppose that $R$ is an $F$-algebra such that $R = [R, R]$. Let $K$
be the algebraic closure of $F$ and $R_K = R \otimes_F K$ the
scalar extension of $R$ to $K$. Then $R_K = [R_K, R_K]$, by
Proposition~\ref{superrings}. Since the property of being a
commutator ring is preserved by homomorphic images, it will
suffice to show that some homomorphic image of $R_K$ is not a
commutator ring; so without loss of generality we may assume that
$R_K$ is a simple $K$-algebra. Then, as a finite-dimensional
simple algebra over an algebraically closed field, $R_K$ is a full
matrix ring. Hence, $[R_K, R_K]$ lies in the kernel of the trace,
by Lemma~\ref{matrix trace}, contradicting $R_K = [R_K, R_K]$.
\end{proof}

From this it follows that, for instance, no PI ring can be a
commutator ring, since every PI ring has an image that is
finite-dimensional over a field.

It is well known that for any $\mathbb Q$-algebra $R$, the $n$th
Weyl algebra $A_n(R)$ satisfies $A_n(R) = [r, A_n(R)]$ for some $r
\in A_n(R)$ (e.g., cf.~\cite{Dixmier}). Actually, it is not hard
to show that writing $A_n(R)= A_1(A_{n-1}(R)) = A_{n-1}(R) \langle
x, y : [x, y] = 1 \rangle$, we can take $r = ax + by +c$ for any
$a, b, c$ in the center $C$ of $R$, such that $aC + bC = C$.
Combining this with another fact about Weyl algebras, we obtain
the following statement.

\begin{proposition}
Let $R$ be a nonzero commutative ring and $n$ a positive integer.
Then the Weyl algebra $A_n(R)$ is a commutator ring if and only if
$R$ is a $\mathbb Q$-algebra.
\end{proposition}

\begin{proof}
If $R$ is not a $\mathbb Q$-algebra, then we may assume that $R$
is a $\mathbb Z /p \mathbb Z$-algebra for some prime $p$, after
passing to a quotient. By a theorem of Revoy in~\cite{Revoy}, for
such a ring $R$, $A_n(R)$ is an Azumaya algebra and hence has a
quotient that is a finite-dimensional algebra over a field.
Therefore, by Proposition~\ref{fd algerbra}, $A_n(R)$ cannot be a
commutator ring.
\end{proof}

We can also prove a more general statement.

\begin{proposition}
Let $n$ be a positive integer, $p$ a prime number, and $R$ a
$\mathbb Z /p \mathbb Z$-algebra. If $R \neq [R, R]$, then $A_n(R)
\neq [A_n(R), A_n(R)]$.
\end{proposition}

\begin{proof}
By induction on $n$, it suffices to prove the proposition for
$A_1(R) = R \langle x, y : [x, y] = 1 \rangle$. Let us denote the
matrix units in $\mathbb M_p(R)$ by $E_{ij}$, and set $X =
\sum_{i=1}^{p-1} E_{i,i+1}$ and $Y = \sum_{i=1}^{p-1} iE_{i+1,i}$.
Then
$$XY - YX = \sum_{i=1}^{p-1} iE_{i,i} - \sum_{i=1}^{p-1}
iE_{i+1,i+1} = \sum_{i=1}^{p} iE_{i,i} - \sum_{i=1}^{p}
(i-1)E_{i,i} = 1.$$
Hence $x \mapsto X$ and $y \mapsto Y$ induces
a ring homomorphism from $A_1(R)$ to $\mathbb M_p(R)$. Now, since
$i^{p-1} \equiv 1 \, (\mathrm{mod} \ p)$ for $i \in \{1, 2, \dots,
p-1\}$, we have
$$\mathrm{trace}((XY)^{p-1}) =
\mathrm{trace}(\sum_{i=1}^{p-1} i^{p-1}E_{i,i}) = \sum_{i=1}^{p-1}
i^{p-1} = p-1.$$
Thus, if $r \in R \backslash [R, R]$, then
$\mathrm{trace}(-r(XY)^{p-1}) = r \in R$, and so, by
Lemma~\ref{matrix trace}, $-r(XY)^{p-1}$ is not a sum of
commutators in $\mathbb M_p(R)$. Consequently, $-r(xy)^{p-1}$ is
not a sum of commutators in $A_1(R)$.
\end{proof}

\section{Endomorphism rings}

We begin with a general result about commutators in matrix rings
and, with its help, provide another way of constructing commutator
rings.

\begin{lemma}
Let $R$ be a ring and $r \in R$ any element. Suppose that $e \in
R$ is an idempotent such that $ere \in [eRe, eRe]_{m_1}$ and $frf
\in [fRf, fRf]_{m_2}$, where $f = 1 - e$ and $m_1, m_2 \in
\mathbb{N}$. Then $r \in [R, R]_{m + 1}$, where $m =
\mathrm{max}(m_1, m_2)$.
\end{lemma}

\begin{proof}
First, we note that for any $x, y, z, w \in R$, $[exe, eye] +
[fzf, fwf] = [exe + fzf, eye + fwf].$ Hence, if $ere \in [eRe,
eRe]_m$ and $frf \in [fRf, fRf]_m$, then $ere + frf \in [R, R]_m$.

Now, write $r = ere + erf + fre + frf$. Then $erf = (-erf)e -
e(-erf)$ and $fre = (fre)e - e(fre)$. Hence, $r = [fre - erf, e] +
ere + frf \in [R, R]_{m + 1}$.
\end{proof}

\begin{proposition}\label{corner}
Let $R$ be a ring, $r \in R$ any element, and $m_1, m_2, \dots,
m_n \in \mathbb{N}$.  Suppose that $e_1, e_2, \dots, e_n \in R$
are orthogonal idempotents such that $1 = e_1 + e_2 + \dots + e_n$
and $e_ire_i \in [e_iRe_i, e_iRe_i]_{m_i}$ for $i = 1, 2, \dots,
n$ $($where each $e_iRe_i$ is viewed as a ring with identity
$e_i$$)$.  Then $r \in [R, R]_{m + n - 1}$, where $m =
\mathrm{max}(m_1, m_2, \dots, m_n)$.
\end{proposition}

\begin{proof}
We will proceed by induction on $n$.  The statement is a tautology
when $n = 1$.  If $n > 1$, let $f = e_1 + e_2 + \dots + e_{n-1}$.
Then $f = f^2$ and $f = 1 - e_n$.  Also, for $i \in \{1, 2, \dots,
n-1 \}$, $$e_i(frf)e_i = e_ire_i \in [e_iRe_i, e_iRe_i]_{m_i} =
[e_i(fRf)e_i, e_i(fRf)e_i]_{m_i},$$ so, by the inductive
hypothesis, $frf \in [fRf, fRf]_{\mathrm{max}(m_1, m_2, \dots,
m_{n-1}) + n - 2}$. But, by assumption, $e_nre_n \in [e_nRe_n,
e_nRe_n]_{m_n}$. Hence, $r \in [R, R]_{m + n - 1}$, by the
preceding lemma and the fact that $\mathrm{max}(\mathrm{max}(m_1,
m_2, \dots, m_{n-1}) + n - 2, m_n) \leq \mathrm{max}(m_1, m_2,
\dots, m_n) + n - 2$.
\end{proof}

\begin{corollary}\label{module sum}
Let $R$ be a ring and $M = M_1 \oplus M_2 \oplus \dots \oplus M_n$
be right R-modules. If each $\mathrm{End}_R(M_i)$ is a commutator
ring, then so is $\mathrm{End}_R(M)$. Also, if for each $i$ there
is some positive integer $m_i$ such that $\mathrm{End}_R(M_i) =
[\mathrm{End}_R(M_i), \mathrm{End}_R(M_i)]_{m_i}$, then
$\mathrm{End}_R(M) = [\mathrm{End}_R(M), \mathrm{End}_R(M)]_{m + n
- 1}$, where $m = \mathrm{max}(m_1, m_2, \dots, m_n)$.
\end{corollary}

Let us now turn to endomorphism rings of infinite direct sums and
products of copies of a fixed module.

\begin{proposition}\label{modules}
Let $R$ be a ring, $N$ a right $R$-module, $\Omega$ an infinite
set, and $M = \bigoplus_{\Omega} N$ or $M = \prod_{\Omega} N$.
Then $\mathrm{End}_R (M) = [x, \mathrm{End}_R (M)]$ for some $x
\in \mathrm{End}_R (M)$. If $\Omega=\mathbb{N}$, then $x$ can be
taken to be the right shift operator.
\end{proposition}

\begin{proof}
Since $\bigoplus_{\Omega} N \cong \bigoplus_{i = 0}^\infty
(\bigoplus_{\Omega} N)$ and $\prod_{\Omega} N \cong \prod_{i =
0}^\infty (\prod_{\Omega} N)$, it suffices to prove the
proposition in the case $\Omega = \mathbb{N}$. Let $x \in
\mathrm{End}_R (M)$ be the right shift operator and $z \in
\mathrm{End}_R (M)$ the left shift operator. Now, consider any
endomorphism $f \in \mathrm{End}_R (M)$, and set $y =
-\sum_{i=0}^\infty x^i fz^{i+1}$. Assuming this summation
converges and using the relation $zx = 1$, we have
$$xy - yx = -\sum_{i=0}^\infty x^{i+1} fz^{i+1} + \sum_{i=0}^\infty x^i
fz^{i+1}x = -\sum_{i=1}^\infty x^i fz^i + \sum_{i=0}^\infty x^i
fz^i = f.$$
It remains to prove convergence of the sum defining
$y$ in the function topology on $\mathrm{End}_R (N)$. (In the case
$N = \bigoplus_{\Omega} M$, it is the topology based on regarding
$N$ as a discrete module, while in the case $N = \prod_{\Omega}
M$, it is the topology constructed using the product topology on
$N$, arising from the discrete topologies on the factors.)

If $M = \bigoplus_{\Omega} N$, then $\bigcup_{i = 1}^\infty
\mathrm{ker}(z^i) = M$. Hence, every element of $M$  is
annihilated by almost all summands of $-\sum_{i=0}^\infty x^i
fz^{i+1}$. If $M = \prod_{\Omega} N$, then given any positive
integer $i$ and any $m \in M$, $x^j fz^{j+1} (m)$ has a nonzero
component in the copy of $N$ indexed by $i$ for only finitely many
values of $j$ (namely, for $j \leq i$). In either case,
$-\sum_{i=0}^\infty x^i fz^{i+1}$ converges in the appropriate
topology on $\mathrm{End}_R (M)$.
\end{proof}

Using similar methods, it can be shown that given any ring $R$,
the ring of infinite matrices over $R$ that are both row-finite
and column-finite is also a commutator ring.

\begin{theorem}
Let $R$ be a ring, $N$ a right $R$-module, $\Omega$ a nonempty
set, and $M = \bigoplus_\Omega N$ or $M = \prod_\Omega N$. Then
$\mathrm{End}_R(M)$ is a commutator ring if and only if either
$\Omega$ is infinite or $\mathrm{End}_R(N)$ is a commutator ring.
\end{theorem}

\begin{proof}
Suppose that $\mathrm{End}_R(M)$ is a commutator ring and $\Omega$
is finite. Then $\mathrm{End}_R(M) \cong \mathbb{M}_n
(\mathrm{End}_R(N))$ for some positive integer $n$. Hence,
$\mathrm{End}_R(N)$ is a commutator ring, by Lemma~\ref{matrix
trace}. The converse follows from the previous proposition if
$\Omega$ is infinite, and from Proposition~\ref{superrings} if
$\Omega$ is finite.
\end{proof}

Incidentally, in the proof of Proposition~\ref{modules}, the only
fact about $M = \bigoplus_{i = 0}^\infty N$ that we used was that
for such a module there are endomorphisms $x, z \in \mathrm{End}_R
(M)$ such that $zx = 1$ and $\bigcup_{i = 1}^\infty
\mathrm{ker}(z^i) = M$. This condition actually characterizes
modules that are infinite direct sums of copies of a module.

\begin{proposition}
Let $R$ be a ring and $M$ a right $R$-module. The following are
equivalent.
\begin{enumerate}
\item[$(1)$] $M = \bigoplus_{i\in\Omega} N_i$ for some infinite
set $\Omega$ and some right $R$-modules $N_i$, such that $N_i
\cong N_j$ for all $i,j \in \Omega$.
\item[$(1')$] $M = \bigoplus_{i\in\Omega} N_i$ for some {\em
countably} infinite set $\Omega$ and some right $R$-modules $N_i$,
such that $N_i \cong N_j$ for all $i,j \in \Omega$.
\item[$(2)$] There exist $x, z \in \mathrm{End}_R (M)$ such that $zx = 1$
and $\bigcup_{i = 1}^\infty \mathrm{ker}(z^i) = M$.
\end{enumerate}
\end{proposition}

\begin{proof}
The equivalence of $(1)$ and $(1')$ is clear. To deduce $(2)$ from
$(1')$, we may assume that $M = \bigoplus_{i = 0}^\infty N_i$,
after well ordering $\Omega$. Then $x$ can be taken to be the
right shift operator and $z$ the left shift operator. To show that
$(2)$ implies $(1')$, let $N_i = x^i \mathrm{ker}(z)$ for each $i
\in \mathbb{N}$. Taking $i \geq 1$ and $a \in \mathrm{ker}(z)$, we
have $zx^i(a) = x^{i-1}(a)$. Hence, left multiplication by $z$
gives a right $R$-module homomorphism $N_i \rightarrow N_{i-1}$,
which is clearly surjective. This homomorphism is also injective,
since $N_i = x(x^{i-1} \mathrm{ker}(z))$ and $zx = 1$. By
induction, $N_i \cong N_j$ for all $i,j \in \mathbb{N}$.

Next, let us show that $\sum_{i = 0}^\infty N_i$ is direct.
Suppose that $a_0 + a_1 + \dots + a_n = 0$, where each $a_i \in
N_i$, $n \geq 1$, and $a_n \neq 0$. For each $i \in \{0, 1, \dots,
n\}$, write $a_i = x^i(b_i)$, where $b_i \in \mathrm{ker}(z)$.
Then $0 = z^n(a_0 + a_1 + \dots + a_n) = z^n(b_0) + z^{n-1}(b_1) +
\dots + z(b_{n-1}) + b_n = b_n$. Hence, $a_n = 0$; a
contradiction.

Finally, we show that given $a \in M$, we have $a \in \bigoplus_{i
= 0}^\infty N_i$. By $(2)$, $a \in \mathrm{ker}(z^n)$ for some
positive integer $n$. If $n = 1$, then $a \in N_0$, so there is
nothing to prove. Otherwise, $za \in \mathrm{ker}(z^{n-1})$, and
we may assume inductively that $za = x^0(b_0) + x^1(b_1) + \dots +
x^{n-1}(b_{n-1})$ for some $b_0, b_1, \dots, b_{n-1} \in
\mathrm{ker}(z)$. Then $a = (a - xza) + xza$, where $a - xza \in
\mathrm{ker}(z) = N_0$, and $xza \in \bigoplus_{i = 1}^n N_i$.
\end{proof}

\section{Traceless matrices}

We now prove our main result about commutators in matrix rings.
This proof uses the same fundamental idea as the one for
Proposition~\ref{modules}.

\begin{theorem}\label{trace}
Let $R$ be a ring and $n$ a positive integer. Then there exist
matrices $X, Y \in \mathbb{M}_n(R)$ such that for all $A \in
\mathbb{M}_n(R)$ having trace $0$, $A \in [X, \mathbb{M}_n(R)] +
[Y, \mathbb{M}_n(R)]$. Specifically, writing $E_{ij}$ for the
matrix units, one can take $X = \sum_{i=1}^{n-1} E_{i+1,i}$ and $Y
= E_{nn}$.
\end{theorem}

\begin{proof}
Write $A = (a_{ij})$, and set $X = \sum_{i=1}^{n-1} E_{i+1,i}$, $Z
= \sum_{i=1}^{n-1} E_{i,i+1}$. Then $ZX = E_{11} + E_{22} + \dots
+ E_{n-1,n-1} = I - E_{nn}$. Hence,
\begin{enumerate}
\item[$(1)$] $ZX = I - E_{nn}.$
\end{enumerate}
Also, for $l \in \{0, 1, \dots, n-1\}$, $$E_{nn}X^lAZ^lE_{nn} =
E_{nn}(\sum_{i=1}^{n-l} E_{i+l,i})A(\sum_{i=1}^{n-l}
E_{i,i+l})E_{nn} \linebreak = E_{n,n-l}AE_{n-l,n} = a_{n-l,
n-l}E_{nn}.$$ Thus,
\begin{enumerate}
\item[$(2)$] $E_{nn}X^lAZ^lE_{nn} = a_{n-l, n-l}E_{nn}.$
\end{enumerate}
Now, let $C = A + XAZ + \dots + X^{n-1}AZ^{n-1}$. Then $[CZ, X] =
CZX - XCZ = C(I - E_{nn}) - XCZ$, by (1). Also, $C - XCZ = A$,
since $X^n = 0$.  Hence, $[CZ, X] = A - CE_{nn}$. We note that
$CE_{nn}$ is an $R$-linear combination of $E_{1n}, E_{2n}, \dots,
E_{n-1, n}$, since $E_{nn}CE_{nn} = (a_{n n} + a_{n-1, n-1} +
\dots + a_{11})E_{nn} = 0$, by (2) and the hypothesis that
$\mathrm{trace}(A) = 0$.

Setting $Y = E_{nn}$, we have for each $i \in \{1, 2, \dots,
n-1\}$, $E_{in} = E_{in}Y - YE_{in}$. Hence $CE_{nn} = [CE_{nn},
Y]$, and therefore $A = [CZ, X] + [CE_{nn}, Y]$.
\end{proof}

\begin{corollary}
Let $R$ be a ring, $n$ a positive integer, $m \in \mathbb{N}$, and
$A \in \mathbb{M}_n(R)$.  If $\mathrm{trace}(A) \in [R, R]_m$,
then $A \in [\mathbb{M}_n(R), \mathbb{M}_n(R)]_{\lceil m/n^2
\rceil + 2}$, where $\lceil m/n^2 \rceil$ denotes the least
integer $\geq m/n^2$.
\end{corollary}

\begin{proof}
Let $r = \mathrm{trace}(A)$. By Lemma~\ref{matrix trace}, there is
a matrix $B \in [\mathbb{M}_n(R), \mathbb{M}_n(R)]_{\lceil m/n^2
\rceil}$ such that $\mathrm{trace}(B) = r$. By the previous
theorem, $A - B$ is a sum of two commutators. Hence, $A \in
[\mathbb{M}_n(R), \mathbb{M}_n(R)]_{\lceil m/n^2 \rceil + 2}$.
\end{proof}

\begin{corollary}\label{Harris-Kap}
Let $R$ be a ring, $n$ a positive integer, and $A \in
\mathbb{M}_n(R)$ a matrix.  Then $A \in [\mathbb{M}_n(R),
\mathbb{M}_n(R)]$ if and only $\mathrm{trace}(A) \in [R, R]$.
\end{corollary}

\begin{proof}
The forward implication was proved in Lemma~\ref{matrix trace},
while the converse follows from the previous corollary.
\end{proof}

This is a generalization of the result in~\cite{Harris} that given
a division ring $D$, a matrix $A \in \mathbb{M}_n(D)$ is a sum of
commutators if and only if its trace is a sum of commutators in
$D$. Actually, Corollary~\ref{Harris-Kap} can also be deduced
rather quickly from Proposition~\ref{corner} by writing $A = B +
C$, where $B \in \mathbb{M}_n(R)$ has zeros on the main diagonal,
$C \in \mathbb{M}_n(R)$ has zeros everywhere off the main
diagonal, and $\mathrm{trace}(C) \in [R, R]$. Then $B \in
[\mathbb{M}_n(R), \mathbb{M}_n(R)]$, by the proposition, and $C$
can be written as a sum of commutators and matrices of the form
$xE_{ii} - xE_{nn} = (xE_{in})E_{ni} - E_{ni}(xE_{in})$.

Let us now give an example of a matrix that has trace $0$ but is
not a commutator, showing that in general the number of
commutators in Theorem~\ref{trace} cannot be decreased to one. A
similar example appears in~\cite{Rosset2}, however, our proof is
considerably shorter. We will require the following lemma in the
process.

\begin{lemma}
Let $F$ be a field and $A, B \in \mathbb{M}_2(F)$ matrices such
that $[A, B] = 0$. Then $\{[A, C] + [B, D] : C, D \in
\mathbb{M}_2(F)\}$ is a subspace of the $F$-vector space
$\mathbb{M}_2(F)$ of dimension at most $2$.
\end{lemma}

\begin{proof}
The result is clear if $A$ and $B$ are both central in
$\mathbb{M}_2(F)$, so we may assume that $A$ is not central. For
all $C \in \mathbb{M}_2 (F)$ let $f(C) = [A, C]$ and $g(C) = [B,
C]$. Then $f$ and $g$ are $F$-vector space endomorphisms of
$\mathbb{M}_2 (F)$. Let $S \subseteq \mathbb{M}_2 (F)$ be the
subspace spanned by the images of $f$ and $g$. We will show that
$S$ is at most $2$-dimensional.

For all $a, b, c \in F$ we have $f(aA + bB + cI) = 0 = g(aA + bB +
cI)$, since $[A, B] = 0$. If $B$ is not in the span of $I$ and
$A$, then $\mathrm{dim}_F (\mathrm{ker}(f)) \geq 3$ and
$\mathrm{dim}_F (\mathrm{ker}(g)) \geq 3$. So $\mathrm{im}(f)$ and
$\mathrm{im}(g)$ are each at most $1$-dimensional, and hence $S$
can be at most $2$-dimensional. Therefore, we may assume that $B =
aA + bI$ for some $a, b \in F$. In this case, $S =
\mathrm{im}(f)$. But, $I, A \in \mathrm{ker}(f)$, and $A$ is not
in the span of $I$, so $\mathrm{im}(f)$ is at most
$2$-dimensional, as desired.
\end{proof}

\begin{proposition}
Let $F$ be a field and $R = F [x_{11}, x_{12}, x_{21}]/(x_{11},
x_{12}, x_{21})^2$. Then the matrix $$X = \left(\begin{array}{rr}
\bar{x}_{11} & \bar{x}_{12} \\
\bar{x}_{21} & -\bar{x}_{11} \\
\end{array}\right) \in \mathbb{M}_2 (R)$$
has trace $0$ but is not a commutator in $\mathbb{M}_2 (R)$.
\end{proposition}

\begin{proof}
Suppose that $X = [A, B]$ for some $A, B \in R$. Viewing $X$, $A$,
and $B$ as polynomials in $\bar{x}_{11}$, $\bar{x}_{12}$, and
$\bar{x}_{21}$, let $A_0, B_0 \in \mathbb{M}_2 (F)$ denote the
degree-$0$ terms of $A$ and $B$, respectively. Then the equation
$X = [A, B]$ tells us that $[A_0, B_0] = 0$ and that the matrices
of the form $[A_0, C] - [B_0, D]$ ($C, D \in \mathbb{M}_2(F)$)
span a $3$-dimensional subspace of $\mathbb{M}_2 (F)$, namely the
subspace of all trace-0 matrices. (For if we denote the
coefficients of $\bar{x}_{ij}$ in $A$ and $B$ by $A_{ij}$ and
$B_{ij}$, respectively, then the coefficient of $\bar{x}_{ij}$ in
$X$ is $[A_0, B_{ij}] + [A_{ij}, B_0]$.) This contradicts the
previous lemma, and so $X \neq [A, B]$ for all $A, B \in R$.
\end{proof}

We can also extend the above result to a more general setting.

\begin{proposition}
Suppose $F$ is a field, $R$ a commutative $F$-algebra, $I
\subseteq R$ an ideal such that $R/I = F$, and $I/I^2$ is at least
$3$-dimensional over $F$. Then there exists a matrix $X \in
\mathbb{M}_2 (R)$ that has trace $0$ but is not a commutator.
\end{proposition}

\begin{proof}
Let $x, y, z \in I$ be such that $\{\bar{x}, \bar{y}, \bar{z}\}$
is $F$-linearly independent in $I/I^2$. Then, there is a
homomorphism $\phi : R \rightarrow F [x_{11}, x_{12},
x_{21}]/(x_{11}, x_{12}, x_{21})^2$ such that $\phi(x) =
\bar{x}_{11}$, $\phi(y) = \bar{x}_{12}$, and $\phi(z) =
\bar{x}_{21}$. The matrix
$$X = \left(\begin{array}{rr}
x & y \\
z & -x \\
\end{array}\right)$$ will then have the desired
properties, since its image in $\mathbb{M}_2 (F [x_{11}, x_{12},
x_{21}]/(x_{11}, x_{12}, x_{21})^2)$ (under the extension of
$\phi$ to a matrix ring homomorphism) is not a commutator, by the
previous proposition.
\end{proof}

\noindent
Department of Mathematics \\
University of California \\
Berkeley, CA 94720 \\
USA \\

\noindent Email: {\tt zak@math.berkeley.edu}

\end{document}